\documentclass[10pt]{article}
\usepackage{amssymb,latexsym,amsmath,graphics,color,std,soul,xcolor}
\usepackage[T1]{fontenc}
\usepackage{amssymb,epsfig}
\usepackage {amsmath}
\usepackage{amsrefs}

\usepackage {bbold}
\usepackage{hyperref}
\usepackage{cleveref}
\usepackage {enumitem}
\usepackage {amsthm}
\usepackage {epsf} 
\usepackage {latexsym}
\newtheorem{theorem}{Theorem}

\theoremstyle{plain}
\newtheorem*{namedthm}{\namedthmname}
\newcounter{namedthm}



\numberwithin{theorem}{section}
\newtheorem{conjecture}[theorem]{Conjecture}
\newtheorem{corollary}[theorem]{Corollary}
\newtheorem{proposition}[theorem]{Proposition}
\newtheorem{lemma}[theorem]{Lemma}

\newtheorem{note}[theorem]{Note}
\newtheorem{definition}[theorem]{Definition}

    \newtheoremstyle{TheoremNum}
        {\topsep}{\topsep}              
        {\itshape}                      
        {}                              
        {\bfseries}                     
        {.}                             
        { }                             
        {\thmname{#1}\thmnote{ \bfseries #3}}
    \theoremstyle{TheoremNum}
    
\newcommand{\paren}[1]{\left(#1\right)}

\begin{document}

\centerline{\Large \bf Proof of the Kresch-Tamvakis Conjecture}  
 
\bigskip

\centerline{\large John S. Caughman\footnote[1]{Corresponding author:  caughman@pdx.edu.} and Taiyo S. Terada}

\bigskip
\bigskip

\centerline{\large \today} 

\bigskip 
\bigskip 

\begin{abstract}
\noindent
In this paper we resolve a conjecture of Kresch and Tamvakis \cite{Kresch}. Our result is the following.

\medskip

\noindent  {\bf Theorem}: For any positive integer $D$ and any integers $i,j$ $(0\leq i,j \leq D)$, the absolute value of the following hypergeometric
series is at most 1:
\begin{equation*} 
 {_4F_3} \left[ \begin{array}{c} -i, \; i+1, \; -j, \; j+1 \\ 1, \; D+2, \; -D \end{array} ; 1 \right].
 \end{equation*}

\noindent
To prove this theorem, we use the Biedenharn-Elliott identity, the theory of Leonard pairs, and the Perron-Frobenius theorem.

\bigskip
\bigskip 



\end{abstract}

\bigskip 
\bigskip

\section{Introduction} 
 
In 2001,  
Kresch and Tamvakis  
conjectured an inequality involving certain terminating ${}_4F_3$ hypergeometric series \cite[Conjecture~2]{Kresch}. In this paper, we prove the conjecture.

\medskip

\noindent
To describe the conjecture, we bring in some notation.
For any real number $a$ and nonnegative integer $n$, define
\begin{equation}\label{poch}
(a)_n = a (a+1) (a+2) \cdots (a+n-1). 
 \end{equation}
Let $z$ denote an indeterminate. Given real numbers $\lbrace a_i \rbrace_{i=1}^4$ and $\lbrace b_i \rbrace_{i=1}^3$, the corresponding ${_4F_3}$ hypergeometric series is defined by
\begin{equation}\label{4f3def}
 {_4F_3} \left[ \begin{array}{c} a_1, \; a_2, \; a_3, \; a_4 \\ b_1, \; b_2, \; b_3 \end{array} ; z \right] = \sum_{n=0}^{\infty} \frac{(a_1)_n (a_2)_n (a_3)_n (a_4)_n}{(b_1)_n (b_2)_n (b_3)_n} \frac{z^n}{n!}.
 \end{equation} 

\noindent
We now state the conjecture of Kresch and Tamvakis. 

\begin{conjecture} \label{KConj} \cite[Conjecture~2]{Kresch}  For any positive integer $D$ and any integers $i,j$ $(0\leq i,j \leq D)$, the absolute value of the following hypergeometric
series is at most 1:
\begin{equation} \label{kresch4F3}
 {_4F_3} \left[ \begin{array}{c} -i, \; i+1, \; -j, \; j+1 \\ 1, \; D+2, \; -D \end{array} ; 1 \right].
 \end{equation}
\end{conjecture}

\begin{note} \rm
    Conjecture~\ref{KConj} is taken from \cite[Conjecture~2]{Kresch} with
\begin{align*}
n=i, \qquad s=j, \qquad T=D+1.
\end{align*}
\end{note} 

\medskip

\noindent
Next we discuss the evidence for Conjecture~\ref{KConj} offered by Kresch, Tamvakis, and others.

\medskip

\noindent
In \cite[Proposition~2]{Kresch},
Kresch and Tamvakis prove that the absolute value of (\ref{kresch4F3}) is at most 1, provided that $i \leq 3$ or $i=D$. 
In \cite[p.~863]{ismail}, Ismail and Simeonov prove that the absolute value of (\ref{kresch4F3}) is at most 1, provided that $i=D-1$ and $D \geq 6$.
They also give asymptotic estimates to further support the conjecture. In \cite{mishev}, Mishev obtains several relations satisfied by the $_4F_3$ hypergeometric series in question. 

\medskip

\noindent
In this paper, we prove Conjecture~\ref{KConj} from scratch, without invoking the above partial results. The following is a statement of our result.

\begin{theorem}\label{thm:main}
    For any positive integer $D$ and any integers $i,j$ $(0\leq i,j \leq D)$, the absolute value of the following hypergeometric
series is at most 1:
\begin{equation*} 
 {_4F_3} \left[ \begin{array}{c} -i, \; i+1, \; -j, \; j+1 \\ 1, \; D+2, \; -D \end{array} ; 1 \right].
 \end{equation*}
\end{theorem}

\noindent
To prove Theorem~\ref{thm:main} we use the following approach. 
For $0 \leq i \leq D$ we define a certain
matrix $B_i \in {\rm Mat}_{D+1}(\mathbb R)$. 
Using the Biedenharn-Elliott identity \cite[p.~356]{biendenharn1985racah}, we show that the entries of $B_i$ are nonnegative.
Using the theory of Leonard pairs \cites{nomura2021idempotent,nomura2021leonard,terwilliger2001two,DefLeonPair2,terwilliger2006algebraic}, we show that the eigenvalues of $B_i$ are $2i+1$ times
\begin{align*}
       {_4F_3} \left[ \begin{array}{c} -i, \; i+1, \; -j, \; j+1 \\ 1, \; D+2, \; -D \end{array} ; 1 \right] \qquad       (0 \leq j \leq D).      
\end{align*}
We also show that the all 1's vector in ${\mathbb{R}}^{D+1}$ is an eigenvector for
$B_i$ with eigenvalue $2i+1$.
Applying the Perron-Frobenius theorem \cite[p.~529]{horn2012matrix}, we show
that the eigenvalues of $B_i$ have absolute value at most $2i+1$.
Using these results, we obtain the proof of Theorem~\ref{thm:main}. 

\medskip

\noindent
This paper is organized as follows. In Section 2, we recall the definition of a Leonard pair and give an example relevant to our work. In Section 3, we use the Leonard pair in Section 2 to
define a sequence of orthogonal polynomials.
In Section 4, we use these orthogonal polynomials to define
the matrices $\lbrace B_i \rbrace_{i=0}^D$. We then compute
the eigenvalues of $\lbrace B_i \rbrace_{i=0}^D$.
In Section 5, we show that the entries of $B_i$ are nonnegative
for $0 \leq i \leq D$.
In Section 6, we use the Perron-Frobenius theorem to prove Theorem~\ref{thm:main}.
In the appendix, we give some details about a key formula in our proof.

\medskip

\noindent
Throughout this paper, the square root of a nonnegative real number is understood to be nonnegative.

\bigskip 

\section{Leonard pairs}

Throughout this paper, $D$ denotes a positive integer. 
Let ${\rm{Mat}}_{D+1}({\mathbb{R}})$ denote the $\mathbb{R}$-algebra of all $(D+1)\times(D+1)$ matrices that have all entries in ${\mathbb{R}}$. We index the rows and columns by $0,1,2,\ldots,D$. Let ${\mathbb{R}}^{D+1}$ denote the vector space over ${\mathbb{R}}$ consisting of $(D+1) \times 1$ matrices that have all entries in ${\mathbb{R}}$. We index the rows by $0,1,2,\ldots,D$.  The algebra ${\rm{Mat}}_{D+1}({\mathbb{R}})$ acts on ${\mathbb{R}}^{D+1}$ by left multiplication.  

\medskip

\noindent
A matrix $B \in {\rm{Mat}}_{D+1}({\mathbb{R}})$ is called {\em tridiagonal} whenever each nonzero entry lies on the diagonal, the subdiagonal, or the superdiagonal. Assume that $B$ is tridiagonal. Then $B$ is called {\em irreducible} whenever each entry on the subdiagonal is nonzero, and each entry on the superdiagonal is nonzero. 

\medskip

\noindent
We now recall the definition of a Leonard pair. Let $V$ denote a vector space over $\mathbb{R}$ with dimension $D+1$.

\begin{definition}\cite{terwilliger2001two} \rm
By a {\it Leonard pair on $V$}, we mean an ordered pair of linear transformations $A:V \rightarrow V$ and $A^*: V \rightarrow V$ that satisfy both {\rm (i), (ii)} below.
\begin{enumerate}
\item  There exists a basis for $V$ with respect to which the matrix representing $A$ is irreducible tridiagonal and the matrix representing $A^*$ is diagonal.
\item  There exists a basis for $V$ with respect to which the matrix representing $A^*$ is irreducible tridiagonal and the matrix representing $A$ is diagonal.
\end{enumerate}
The above Leonard pair $A,A^*$ is said to be {\it over $\mathbb{R}$}.
\end{definition}

\begin{note} \rm
    According to a common notational convention, $A^*$ denotes the conjugate-transpose of $A$. We are not using this convention. In a Leonard pair $A, A^*$ the linear transformations $A$ and $A^*$ are arbitrary subject to (i), (ii) above. 
\end{note} 

\noindent 
Our next goal is to give an example of a Leonard pair. To do so, we give two definitions.

\smallskip

\begin{definition}\label{maindef} \rm
Define
\begin{align}\label{abcDef}
c_i &= \frac{3(D-i+1)i(D+i+1)}{D(D+2)(2i+1)} &(1 \leq i \leq D),\\
a_i &= \frac{3i(i+1)}{D(D+2)}  &(0 \leq i \leq D),\\
b_i &= \frac{3(D-i)(i+1)(D+i+2)}{D(D+2)(2i+1)} &(0 \leq i \leq D-1),\\
\theta_i &= 3-2a_i & (0 \leq i \leq D). 
\end{align}
We remark that the scalars $\{ \theta_i \}_{i=0}^D$ are mutually distinct. \\

\noindent
Let $A, A^*$ denote the following matrices in ${\rm{Mat}}_{D+1}({\mathbb{R}})$:
\begin{equation}\label{def:KTpair} 
A = \left( \begin{array}{ccccc}
a_0 & b_0 &        &        & {\bf{0}} \\
c_1 & a_1 & b_1    &        &  \\
    & \ddots & \ddots & \ddots &    \\
    &         & c_{D-1}& a_{D-1}&  b_{D-1}\\
{\bf{0}} &        &        & c_{D}& a_{D}
\end{array}
\right),
\hspace{3em}
A^* = \left( \begin{array}{ccccc}
\theta_0 & &        &        & {\bf{0}} \\
 & \theta_1 &   &        &  \\
    &  & \ddots &  &    \\
    &         & & \theta_{D-1}&  \\
{\bf{0}} &        &        & & \theta_{D}
\end{array}
\right). 
\end{equation} 
\end{definition} 

\bigskip

\begin{definition}\label{def:P} \rm
We define a matrix $P \in {\rm{Mat}}_{D+1}({\mathbb{R}})$ with the following entries:
    \begin{equation}\label{Pdef}
P_{i,j} =(2j+1)\, {_4F_3}\left[ \begin{array}{c} -i, \; i+1, \; -j, \; j+1 \\ 1, \; D+2, \; -D \end{array} ; 1 \right]  \qquad (0 \leq i,j \leq D).
\end{equation}
\end{definition}

\medskip

\begin{lemma}\label{lem:A-As-form-leonard-pair} {\rm (\cite[Ex.~5.10]{DefLeonPair2} and
\cite[Thm.~4.9]{terwilliger2006algebraic})}
The following hold:
\begin{enumerate}
    \item $P^2=(D+1)^2I$;
    \item $PA=A^*P$;
    \item $PA^*=AP$;
    \item the pair $A,A^*$ is a Leonard pair over $\mathbb{R}$.
\end{enumerate}
\end{lemma}
\begin{proof}
    Calculations {\em (i)}--{\em(iii)} are the following special case of \cite[Ex.~5.10]{DefLeonPair2} and \cite[Thm.~4.9]{terwilliger2006algebraic}: $$d=D, \qquad \theta_0=\theta_0^*=3, \qquad s=s^*=r_1=0, \qquad r_2=D+1, \qquad h=h^*=\frac{-6}{D(D+2)}.$$
    Item {\em(iv)} follows from items {\em (i)}--{\em(iii)}.
\end{proof}

\noindent
The Leonard pairs from \cite[Ex.~5.10]{DefLeonPair2} are said to have Racah type. So the Leonard pair $A,A^*$ in Lemma~\ref{lem:A-As-form-leonard-pair} has Racah type. This Leonard pair is self-dual in the sense of \cite[p.~5]{nomura2021leonard}.

\bigskip 

\section{Some orthogonal polynomials}

\noindent 
In this section we interpret Conjecture~\ref{KConj} in terms of orthogonal polynomials. \\


\noindent
Let $\lambda$ denote an indeterminate. Let ${\mathbb R}[\lambda]$ denote the $\mathbb{R}$-algebra of polynomials in $\lambda$ that have all coefficients in ${\mathbb{R}}$.

\begin{definition}\label{def:ui} \rm With reference to Definition \ref{maindef}, let
$u_0(\lambda), u_1(\lambda), \ldots, u_D(\lambda)$
denote the polynomials in ${\mathbb R}[\lambda]$ that satisfy:
$$ u_0(\lambda)=1,\qquad u_1(\lambda)=\lambda/3,$$
\begin{equation}\label{urecur}
    \lambda u_i(\lambda) = b_i u_{i+1}(\lambda) + a_i u_i(\lambda) + c_i u_{i-1}(\lambda) \hspace{2em} (1 \leq i \leq D-1).
\end{equation}  
\end{definition}

\medskip

\noindent
Note that the polynomial $u_i(\lambda)$ has degree exactly $i$ for $0\leq i \leq D$. 

\medskip

\noindent
By \cite[Ex.~5.10]{DefLeonPair2}, the polynomials $\lbrace u_i(\lambda)\rbrace_{i=0}^D$ are a special case of the Racah polynomials. Also by \cite[Ex.~5.10]{DefLeonPair2},
\begin{equation}\label{uijF}
 u_i(\theta_j) = {_4F_3} \left[ \begin{array}{c} -i, \; i+1, \; -j, \; j+1 \\ 1, \; D+2, \; -D \end{array} ; 1 \right]  \hspace{3em} (0 \leq i,j \leq D).
\end{equation} 

\bigskip

\begin{lemma} The following hold:
\begin{enumerate}
    \item $u_i(\theta_j) = u_j(\theta_i) \qquad (0 \leq i,j \leq D)$;
    \item $u_i(\theta_0) = 1 \qquad (0 \leq i \leq D)$;
    \item $u_0(\theta_j) = 1 \qquad (0 \leq j \leq D)$.
\end{enumerate}  
\end{lemma} 
\begin{proof}
Each of $(i)$--$(iii)$ is immediate from (\ref{uijF}).    
\end{proof}

\noindent
In light of Equation~(\ref{uijF}), Conjecture~\ref{KConj} asserts that
\begin{equation}\label{conj2}
 |u_i(\theta_j)| \leq 1 \hspace{3em} (0 \leq i,j \leq D).
\end{equation}

\noindent
To prove (\ref{conj2}) it will be useful to adjust the normalization of the polynomials $u_i(\lambda)$. 

\medskip

\noindent
Define 
\begin{equation}\label{def:ki}
    k_i=\frac{b_0b_1\cdots b_{i-1}}{c_1c_2\cdots c_i} \qquad \qquad (0 \leq i \leq D).
\end{equation}
One checks that \begin{equation}\label{kidef}
 k_i = 2i+1 \hspace{2em} (0 \leq i \leq D).
 \end{equation}

\noindent

\begin{definition}\label{def:vi} \rm With reference to Definition~\ref{def:ui}, let
    \begin{equation}\label{uvreln}
 v_i(\lambda) = k_i u_i(\lambda)  \hspace{2em} (0 \leq i \leq D).
 \end{equation}
\end{definition}

\noindent
By construction,
\begin{equation}\label{uvrelnTH}
 v_i(\theta_j) = k_i u_i(\theta_j)  \hspace{2em} (0 \leq i,j \leq D).
 \end{equation}

\noindent
The polynomials $v_i(\lambda)$ satisfy the following three-term recurrence.
\begin{lemma}\cite[Lem.~3.11]{terwilliger2006algebraic}\label{lem:vi}   We have
$$v_0(\lambda)=1,\qquad v_1(\lambda)=\lambda,$$
\begin{equation}\label{vrecur}
 \lambda v_i(\lambda) = c_{i+1} v_{i+1}(\lambda) + a_i v_i(\lambda) + b_{i-1} v_{i-1}(\lambda) \hspace{2em} (1 \leq i \leq D-1).
 \end{equation}
\end{lemma}

\begin{lemma}\label{lem:P-in-terms-of-vi}For $0\leq i,j \leq D$ we have
\begin{equation} \label{eq:Pv}
    P_{i,j} = v_j(\theta_i).
\end{equation}
\end{lemma}
\begin{proof}
    Immediate by (\ref{Pdef}),(\ref{uijF}),(\ref{kidef}), and (\ref{uvrelnTH}).
\end{proof}

\medskip


\noindent
We emphasize two special cases of (\ref{eq:Pv}).

\begin{lemma}\label{lem:P-first-row-col} The following hold:
\begin{enumerate}
    \item $P_{i,0} = 1 \qquad \qquad (0 \leq i \leq D);$
    \item $P_{0,j} = k_j \qquad \qquad (0 \leq j \leq D).$
\end{enumerate}
\end{lemma}
\begin{proof}
Immediate from (\ref{uvrelnTH}) and (\ref{eq:Pv}).
\end{proof}
 

\noindent
We have some comments about the parameters (\ref{def:ki}). For notational convenience, define
\begin{equation}\label{eq:nu}
 \nu= (D+1)^2.
 \end{equation}
 By (\ref{kidef}),  $$\sum_{i=0}^Dk_i = \nu.$$


\noindent
Next, we state the orthogonality relations for the polynomials $\{u_i(\lambda)\}_{i=0}^D$.

\begin{lemma}
\cite[p.~282]{terwilliger2006algebraic} For integers $0 \leq n,m \leq D$ we have
\begin{align}
\label{eq:orthogonality_u-self-dual}
		\sum_{j=0}^D k_j u_n(\theta_j)u_m(\theta_j)&=\nu k_n^{-1}\delta_{n,m};\\
\label{eq:orthogonality2_u-self-dual}
		\sum_{j=0}^Dk_j u_j(\theta_n)u_j(\theta_m)&=\nu k_n^{-1}\delta_{n,m}.
\end{align}
\end{lemma}

\medskip

\noindent
Next, we state the orthogonality relations for the polynomials $\{v_i(\lambda)\}_{i=0}^D$.

\begin{lemma}
\cite[p.~281]{terwilliger2006algebraic} For integers $0 \leq n,m \leq D$ we have
\begin{align}
\label{eq:orthogonality_v-self-dual}
		\sum_{j=0}^D k_j v_n(\theta_j)v_m(\theta_j)&=\nu k_n\delta_{n,m}; \\
\label{eq:orthogonality2_v-self-dual}
	\sum_{j=0}^D k_j^{-1} v_j(\theta_n)v_j(\theta_m)&=\nu k_n^{-1}\delta_{n,m}.
\end{align}
\end{lemma}

\bigskip

\section{Two commutative subalgebras of ${\rm Mat}_{D+1}({\mathbb R})$}

\noindent
We continue to discuss the Leonard pair $A, A^*$ from Definition~\ref{maindef}.

\noindent 
\begin{definition}\label{def:M-Ms-alg} \rm
    Let $M$ denote the subalgebra of ${\rm Mat}_{D+1}({\mathbb R})$ generated by $A$. Let $M^*$ denote the subalgebra of ${\rm Mat}_{D+1}({\mathbb R})$ generated by $A^*$.
\end{definition}

\noindent
In this section, we describe a basis for $M$ and a basis for $M^*$. 

\begin{definition}\label{def:Bi} \rm For $0\leq i \leq D$ define
$$B_i = v_i(A), \qquad\qquad  B_i^*=v_i(A^*),$$
where $v_i(\lambda)$ is from (\ref{uvreln}).
 
\end{definition}


\begin{lemma}\label{lem:PBi-BisP} For $0\leq i \leq D$ we have
$$PB_i=B_i^*P, \qquad \qquad PB_i^*=B_iP.$$
\end{lemma}
\begin{proof}
    By Lemma~\ref{lem:A-As-form-leonard-pair}, Definition~\ref{def:Bi}, and linear algebra.
\end{proof}

\noindent
Lemma~\ref{lem:PBi-BisP} tells us that for integers $0\leq i,j \leq D$, column $j$ of $P$ is an eigenvector of $B_i$ with eigenvalue $v_i(\theta_j)$.
We emphasize one special case. Let $\mathbb{1}$ denote the vector in ${\mathbb{R}}^{D+1}$ that has all entries 1.

\begin{lemma}\label{lem:bi-eig-ki}
For $0 \leq i \leq D$ the vector $\mathbb{1}$ is an eigenvector for $B_i$ with eigenvalue $k_i$. 
\end{lemma} 
\begin{proof}
Immediate from Lemma~\ref{lem:P-first-row-col} and Lemma~\ref{lem:PBi-BisP}.    
\end{proof}

\begin{lemma}
    The matrices $\{B_i\}_{i=0}^D$ form a basis for $M$. The matrices $\{B_i^*\}_{i=0}^D$ form a basis for $M^*$.
\end{lemma}
\begin{proof}
    By Definition \ref{maindef}, the matrix $A^*$ has $D+1$ distinct eigenvalues, so $M^*$ has dimension $D+1$. By Definition \ref{def:Bi}, the matrices $\{B_i^*\}_{i=0}^D$ belong to $M^*$. By these comments, the matrices $\{B_i^*\}_{i=0}^D$ form a basis for $M^*$. We have now verified the second assertion. The first assertion follows from this and Lemma~\ref{lem:PBi-BisP}.
\end{proof}

\noindent
Next we discuss the entries of the matrices $\{ B_i \}_{i=0}^D$. The following definition will be convenient.

\begin{definition}\label{def:phijDef} \rm
    For $0 \leq h,i,j \leq D$ let $p^h_{i,j}$ denote the $(h,j)$-entry of $B_i$. In other words, 
\begin{equation}
\label{phijDef}
p^h_{i,j} = (B_i)_{h,j}.
\end{equation}
\end{definition}


\noindent
 We have a comment about the scalars $p^h_{i,j}$. 
\begin{lemma}\label{lem:struct} \cite[Lem.~4.19]{nomura2021leonard} For $0\leq i,j \leq D$ we have
    \begin{equation}\label{eq:struct}
B_i B_j = \sum_{h=0}^D p^h_{i,j} B_h, \qquad \qquad B^*_i B^*_j = \sum_{h=0}^D p^h_{i,j} B^*_h.
\end{equation}
\end{lemma}
 
\noindent
The scalars $p^h_{i,j}$ can be computed using the following result. This result is from \cite{nomura2021idempotent}; we include a proof for the sake of completeness.

\begin{proposition}\label{prop:phij-ui-triple-product}\cite[Lem.~12.12]{nomura2021idempotent} 
For $0\leq h,i,j \leq D$ we have
    \begin{equation}\label{eq:verlinde}
         p^h_{i,j}=\frac{k_ik_j}{\nu}\sum\limits_{t=0}^{D}k_t u_t(\theta_i)u_t(\theta_j)u_t(\theta_h).
     \end{equation}
\end{proposition}
\begin{proof}
   We invoke Equation (\ref{phijDef}).  By Lemma~\ref{lem:A-As-form-leonard-pair}(i) and Lemma~\ref{lem:PBi-BisP} we have that $B_i=\nu^{-1}PB_i^*P$. Recall that the matrix $P$ has entries $P_{i,j}=k_j u_j(\theta_i)$. We also have $B^*_i = v_i(A^*)$ and $A^*=\text{diag}(\theta_0,\theta_1,\ldots,\theta_D)$.  
    Evaluating (\ref{phijDef}) using these comments, we obtain the result. 
\end{proof}

\medskip

\noindent
We have a comment about Proposition~\ref{prop:phij-ui-triple-product}.

\begin{lemma}\label{lem:phij-sym} For $0\leq h,i,j \leq D$ we have
\begin{equation}\label{eq:phij-symmetry}
 p^h_{i,j}= p^h_{j,i}, \qquad\qquad k_hp^h_{i,j}=k_jp^j_{h,i}=k_ip^i_{j,h}.
\end{equation}
\end{lemma}
\begin{proof}
    Immediate from (\ref{eq:verlinde}).
\end{proof}


\bigskip 

\section{The nonnegativity of the $p^h_{i,j}$}

Our next goal is to show that $p^h_{i,j} \geq0$ for $0\leq h,i,j\leq D$. To obtain this inequality, we use the Biedenharn-Elliott identity \cite[p.~356]{biendenharn1985racah}.

\medskip 

\noindent
Recall the natural numbers $\mathbb{N}=\{0,1,2,3,\ldots\}$. Note that 
$\frac{1}{2}\mathbb{N}=\{0,\frac{1}{2},1,\frac{3}{2},2,\frac{5}{2},\ldots\}$.

\begin{definition}\label{def:admissible} \rm
Given $a,b,c \in \frac{1}{2}\mathbb{N}$, we say that the triple $(a,b,c)$ is {\it admissible} whenever 
$a+b+c\in \mathbb{N}$ and
\begin{equation}\label{eq:tri-ineq}
    a\leq b+c,\qquad \qquad b \leq c+a, \qquad \qquad c\leq a+b.
\end{equation}
\end{definition}

\noindent
\begin{definition}\label{def:tri-coeff} \rm
Referring to Definition~\ref{def:admissible}, assume that $(a,b,c)$ is admissible. Define
 \begin{equation}\label{eq:triangle-coeff}
     \Delta(a,b,c)=\left(\frac{(a+b-c)!(b+c-a)!(c+a-b)!}{(a+b+c+1)!}\right)^{\frac{1}{2}}.
 \end{equation} 
\end{definition}

\noindent
Next, we recall the Racah coefficients.
 
 \begin{definition}\label{Wdef} \rm
 (\cite[Eq.~5.11.4]{biendenharn1985racah} and \cite[p.~1063]{messiah1962racah})  For $a,b,c,d,e,f\in\frac{1}{2}\mathbb{N}$, we define a real number $W(a,b,c,d;e,f)$ as follows.

\medskip

\noindent
First assume that each of $(a,b,e)$, $(c,d,e)$, $(a,c,f)$, $(b,d,f)$ is admissible. Then
 \begin{align}
 \begin{split}\label{eq:bidenharn-louck-racah-coeff-def}      W(a,b,c,d;e,f)&=\frac{\Delta(a,b,e)\Delta(c,d,e)\Delta(a,c,f)\Delta(b,d,f) (\beta_1 +1)!(-1)^{\beta_1-(a+b+c+d)}}{(\beta_2-\beta_1)!(\beta_3-\beta_1)!(\beta_1-\alpha_1)!(\beta_1-\alpha_2)!(\beta_1-\alpha_3)!(\beta_1-\alpha_4)!}\\
    &\quad\times{}_4F_3\left[\begin{matrix}
						\alpha_1-\beta_1 ,\, \alpha_2-\beta_1 ,\, \alpha_3-\beta_1 ,\, \alpha_4-\beta_1\\
						-\beta_1-1 ,\, \beta_2-\beta_1+1 ,\, \beta_3-\beta_1+1 &
					\end{matrix};\,1\right],
\end{split}
 \end{align}
 where 
 $$(\alpha_1,\alpha_2,\alpha_3,\alpha_4) = \mbox{any permutation of } \, (a+b+e, \; c+d+e, \; a+c+f, \; b+d+f),$$
 and where 
 $$\beta_1=\min(a+b+c+d, \; a+d+e+f, \; b+c+e+f),$$ and $\beta_2, \beta_3$ are the other two values in the triple $(a+b+c+d,a+d+e+f,b+c+e+f)$ in either order. 

\medskip

 \noindent
 Next assume that $(a,b,e)$, $(c,d,e)$, $(a,c,f)$, $(b,d,f)$, are not all admissible. Then
 \begin{equation}\label{eq:tri-racah}
     W(a,b,c,d;e,f)=0.
 \end{equation} 
 We call $W(a,b,c,d;e,f)$ the {\it Racah coefficient} associated with $a,b,c,d,e,f$.
\end{definition}


    \medskip

\noindent
Let $0\leq h,i,j\leq D$. In order to show that $p^h_{i,j} \geq 0$, we will show that
    \begin{align*}
        p^h_{i,j}=(2i+1)(2j+1)(D+1) \Bigl( W\left({\textstyle \frac{D}{2},\frac{D}{2},i,h;j,\frac{D}{2}}\right)\Bigr)^2.
    \end{align*}
    We will use the Biedenharn-Elliott identity.

\begin{proposition}\label{prop:biedenharn-elliott} {\rm (Biedenharn-Elliott identity \cite[p.~356]{biendenharn1985racah})}  Let $a,a',b,b',c,c',e,f,g\in\frac{1}{2}\mathbb{N}$. Then 
\begin{equation}
    \begin{split}\label{eq:biedenharn-elliott} 
    \sum_{d\in \frac{1}{2}\mathbb{N}}(-1)^{c+c'-d}(2d+1)&W(b,b',c,c';d,e)W(a,a',c,c';d,f)W(a,a',b,b';d,g) \\
    &=(-1)^{e+f-g}W(a,b,f,e;g,c)W(a',b',f,e;g,c').
\end{split}
\end{equation}

\end{proposition}

 
\noindent
In order to evaluate the Racah coefficients in the Biedenharn-Elliott identity, we will use the following transformation formula of Whipple.

\begin{proposition}\label{prop:whipple} {\rm (Whipple transformation \cite[p.~49]{gasper2004basic})} For integers $p,q,a_1,a_2,r,b_1,b_2$ we have 
\begin{equation}\label{eq:whipple}
    {}_4F_3\left[\begin{matrix}
						-p ,\, q ,\, a_1 ,\, a_2\\
						r ,\, b_1 ,\, b_2 &
					\end{matrix};\,1\right]=\frac{(b_1-q)_p(b_2-q)_p}{(b_1)_p(b_2)_p}{}_4F_3\left[\begin{matrix}
						-p ,\, q ,\, r-a_1 ,\, r-a_2\\
						r ,\, 1+q-b_1-p ,\, 1+q-b_2-p &
					\end{matrix};\,1\right],
\end{equation}
provided that $p\geq0$ and $q+a_1+a_2+1=r+b_1+b_2+p$.
\end{proposition} 

\medskip

\noindent
We are interested in the following Racah coefficient. For $0\leq i,j\leq D$ consider $$W\left({\textstyle \frac{D}{2},\frac{D}{2},\frac{D}{2},\frac{D}{2};i,j}\right).$$ Evaluating this Racah coefficient using Definition~\ref{Wdef} we get a scalar multiple of a certain ${}_4 F_3$ hypergeometric series. Applying several Whipple transformations to this hypergeometric series, we get the following result as we will see.

 \begin{proposition}\label{prop:racah-coeff-kresch-case}
 For integers $0\leq i,j\leq D$ we have
    \begin{equation}\label{eq:racah-coeff-kresch-case}
    W\left({\textstyle \frac{D}{2},\frac{D}{2},\frac{D}{2},\frac{D}{2};i,j}\right)=\frac{(-1)^{i+j-D}}{D+1}{}_4F_3\left[\begin{matrix}
						-i ,\, i+1 ,\, -j ,\, j+1\\
						1 ,\, D+2 ,\, -D &
					\end{matrix};\,1\right].
    \end{equation}
 \end{proposition}

\begin{proof} 
To evaluate $W\left({\textstyle \frac{D}{2},\frac{D}{2},\frac{D}{2},\frac{D}{2};i,j}\right)$, we will consider two cases: $i+j\leq D$ and $i+j>D$.\\

\noindent
{\bf Case $i+j\leq D$}. In this case, from (\ref{eq:bidenharn-louck-racah-coeff-def}) we get $\beta_1=D+i+j$, $\beta_2=2D$, $\beta_3=D+i+j$, $\alpha_1=\alpha_2=D+i$, $\alpha_3=\alpha_4=D+j$. The hypergeometric term in (\ref{eq:bidenharn-louck-racah-coeff-def}), after rearranging the upper indices, becomes
\begin{equation}\label{case1forFversion1}
    {}_4F_3\left[\begin{matrix}
						-i ,\, -i ,\, -j ,\, -j\\
						-D-i-j-1 ,\, D-i-j+1 ,\, 1 &
					\end{matrix};\,1\right].
\end{equation}
The coefficient in (\ref{eq:bidenharn-louck-racah-coeff-def}) is 
\begin{align}\label{eq:racah-coeff-n+j-leq-D-0}
&\frac{\Bigl(\Delta\left({\textstyle\frac{D}{2},\frac{D}{2},i}\right)\Bigr)^2\Bigl(\Delta\left({\textstyle\frac{D}{2},\frac{D}{2},j}\right)\Bigr)^2(D+i+j+1)!(-1)^{i+j-D}}{(D-i-j)!(j!)^2(i!)^2}\nonumber\\
&\qquad \qquad=\frac{(D-i)!(i!)^2(D-j)!(j!)^2(D+i+j+1)!(-1)^{i+j-D}}{(D+i+1)!(D+j+1)!(D-i-j)!(j!)^2(i!)^2}.
\end{align}
The expression (\ref{eq:racah-coeff-n+j-leq-D-0}) is equal to
\begin{equation}\label{eq:racah-coeff-n+j-leq-D}
\frac{(D-i)!(D-j)!(D+i+j+1)!(-1)^{i+j-D}}{(D+i+1)!(D+j+1)!(D-i-j)!}.
\end{equation}
Performing a Whipple transformation (\ref{eq:whipple}) with the substitutions
$-p=-i$,  $q=-j$, $a_1=-i$, $a_2=-j$, $r=1$, $b_1=-D-i-j-1$, $b_2=D-i-j+1$, the hypergeometric component in (\ref{case1forFversion1}), after rearranging lower indices, becomes
\begin{equation}\label{case1forFversion2}
    {}_4F_3\left[\begin{matrix}
						-i ,\, i+1 ,\, -j ,\, j+1\\
						1 ,\, D+2 ,\, -D &
					\end{matrix};\,1\right].
\end{equation}
The coefficient contribution from the Whipple transformation is
\begin{equation}\label{eq:whipple-coeff-n+j-leq-D}
    \frac{(-D-i-1)_i(D-i+1)_i}{(-D-i-j-1)_i(D-i-j+1)_i}=\frac{(-1)^i(D+i+1)!}{(D+1)!}\frac{D!}{(D-i)!}\frac{(D+j+1)!}{(-1)^i(D+i+j+1)!}\frac{(D-i-j)!}{(D-j)!}.
\end{equation}
We see that coefficients (\ref{eq:racah-coeff-n+j-leq-D}) and (\ref{eq:whipple-coeff-n+j-leq-D}) multiply to $\frac{(-1)^{i+j-D}}{D+1}$, as desired.\\

\noindent
{\bf Case $i+j>D$}.
 In this case, from (\ref{eq:bidenharn-louck-racah-coeff-def}) we get $\beta_1=2D$, $\beta_2=D+i+j$, $\beta_3=D+i+j$, $\alpha_1=\alpha_2=D+i$, $\alpha_3=\alpha_4=D+j$. The hypergeometric term in (\ref{eq:bidenharn-louck-racah-coeff-def}) becomes
\begin{equation}\label{case2forFversion1}
    {}_4F_3\left[\begin{matrix}
						i-D ,\, i-D ,\, j-D ,\, j-D\\
						-2D-1 ,\, i+j-D+1 ,\, i+j-D+1 &
					\end{matrix};\,1\right].
\end{equation}
The coefficient in (\ref{eq:bidenharn-louck-racah-coeff-def}) is 
\begin{equation}\label{eq:racah-coeff-n+j-geq-D-0}
    \frac{\Bigl(\Delta\left({\textstyle\frac{D}{2},\frac{D}{2},i}\right)\Bigr)^2\Bigl(\Delta\left({\textstyle\frac{D}{2},\frac{D}{2},j}\right)\Bigr)^2(2D+1)!}{\bigl((i+j-D)!\bigr)^2\bigl((D-i)!\bigr)^2\bigl((D-j)!\bigr)^2}=\frac{(D-i)!(i!)^2(D-j)!(j!)^2(2D+1)!}{(D+i+1)!(D+j+1)!\bigl((i+j-D)!(D-i)!(D-j)!\bigr)^2}.
\end{equation}
The expression (\ref{eq:racah-coeff-n+j-geq-D-0}) is equal to
\begin{equation}\label{eq:racah-coeff-n+j-geq-D}
C_0=\frac{(i!)^2(j!)^2(2D+1)!}{(D+i+1)!(D+j+1)!\bigl((i+j-D)!\bigr)^2(D-i)!(D-j)!}.
\end{equation}
Now we will perform three Whipple transformations. For each one we list the indices chosen $-p$, $q$, $a_1$, $a_2$, $r$, $b_1$, $b_2$, the resulting hypergeometric term (with possible rearranging of some upper indices), and the coefficient contribution, $C_i$, from the corresponding Whipple transformation. 

\begin{enumerate}[label=\arabic*.]
    \item Using $-p=i-D$, $q=j-D$,  $a_1=i-D$, $a_2=j-D$, $r=i+j-D+1$, $b_1=-2D-1$, $b_2=i+j-D+1$:

\begin{equation}\label{case2forFversion2}
    {}_4F_3\left[\begin{matrix}
						i-D ,\, i+1 ,\, j-D ,\, j+1\\
						i+j+2 ,\, -D ,\, i+j-D+1 &
					\end{matrix};\,1\right],
\end{equation}

\begin{align}
    C_1&=\frac{(-D-j-1)_{D-i}(i+1)_{D-i}}{(-2D-1)_{D-i}(i+j-D+1)_{D-i}}\nonumber\\
    &=\frac{(-1)^{D-i}(D+j+1)!}{(i+j+1)!}\frac{D!}{i!}\frac{(D+i+1)!}{(-1)^{D-i}(2D+1)!}\frac{(i+j-D)!}{j!}
    .\label{eq:whipple-coeff-i+j-geq-D-1}
\end{align}

\item Using $-p=i-D$, $q=j+1$, $a_1=i+1$, $a_2=j-D$, $r=-D$, $b_1=i+j+2$, $b_2=i+j-D+1$:

\begin{equation}\label{case2forFversion3}
    {}_4F_3\left[\begin{matrix}
						i-D ,\, -D-i-1 ,\, -j ,\, j+1\\
						 -D ,\, -D ,\, 1 &
					\end{matrix};\,1\right],
\end{equation}

\begin{align}
    C_2&=\frac{(i+1)_{D-i}(i-D)_{D-i}}{(i+j+2)_{D-i}(i+j-D+1)_{D-i}}\nonumber\\
    &=\frac{D!}{i!}(-1)^{D-i}(D-i)!\frac{(i+j+1)!}{(D+j+1)!}\frac{(i+j-D)!}{j!}.\label{eq:whipple-coeff-i+j-geq-D-2}
\end{align}

\item Using $-p=-j$, $q=j+1$, $a_1=i-D$, $a_2=-D-i-1$, $r=-D$, $b_1=-D$, $b_2=1$:

\begin{equation}\label{case2forFversion4}
    {}_4F_3\left[\begin{matrix}
						-i ,\, i+1 ,\, -j ,\, j+1\\
						-D ,\, D+2 ,\, 1 &
					\end{matrix};\,1\right]
     ={}_4F_3\left[\begin{matrix}
						-i ,\, i+1 ,\, -j ,\, j+1\\
						1 ,\, D+2 ,\, -D &
					\end{matrix};\,1\right],
\end{equation}

\begin{align}
    C_3&=\frac{(-D-j-1)_{j}(-j)_{j}}{(-D)_{j}(1)_{j}}\nonumber\\
    &=\frac{(-1)^{j}(D+j+1)!}{(D+1)!}(-1)^j j!\frac{(D-j)!}{(-1)^{j}D!}\frac{1}{j!}
    .\label{eq:whipple-coeff-i+j-geq-D-3}
\end{align}
\end{enumerate}
Combining coefficients we see that $C_0C_1C_2C_3=\frac{(-1)^{D-i+j}}{D+1}=\frac{(-1)^{i+j-D}}{D+1}$, since $i,j,D$ are integers.
\end{proof}

\noindent
We now evaluate the Biedenharn-Elliott identity using Proposition~\ref{prop:racah-coeff-kresch-case}.

\begin{proposition}\label{prop:biedenharn-elliott-kresch}  For integers $0\leq h,i,j\leq D$ we have 
\begin{equation}\label{eq:biedenharn-elliott-kresch}
    \sum_{t=0}^D(2t+1)u_t(\theta_h)u_t(\theta_i)u_t(\theta_j)=(D+1)^3\Bigl(W\left({\textstyle \frac{D}{2},\frac{D}{2},i,h;j,\frac{D}{2}}\right)\Bigr)^2.
\end{equation}
\end{proposition}

\begin{proof}
    First we apply Proposition \ref{prop:biedenharn-elliott} with $a=a'=b=b'=c=c'=\frac{D}{2}$, $e=h$, $f=i$, $g=j$, and $d=t$ to obtain
\begin{align}\label{eq:using-BE} 
    \sum_{t\in \frac{1}{2}\mathbb{N}}(-1)^{D-t}(2t+1)&W({\textstyle \frac{D}{2},\frac{D}{2},\frac{D}{2},\frac{D}{2};t,h})W({\textstyle \frac{D}{2},\frac{D}{2},\frac{D}{2},\frac{D}{2};t,i})W({\textstyle \frac{D}{2},\frac{D}{2},\frac{D}{2},\frac{D}{2};t,j})\nonumber\\
    &=(-1)^{h+i-j}W\left({\textstyle \frac{D}{2},\frac{D}{2},i,h;j,\frac{D}{2}}\right)W\left({\textstyle \frac{D}{2},\frac{D}{2},i,h;j,\frac{D}{2}}\right).
\end{align}
    Note that $\frac{D}{2}+\frac{D}{2}+t$ is an integer if and only if $t$ is an integer. So by (\ref{eq:tri-racah}), 
    the terms of the sum vanish in which $t$ is not an integer or $t>D$.
    By Proposition~\ref{prop:racah-coeff-kresch-case} and (\ref{uijF}), the left hand side of (\ref{eq:using-BE}) becomes 
    \begin{align}
        &\sum_{t=0}^D(-1)^{D-t}(2t+1)\frac{(-1)^{t+h-D}u_t(\theta_h)}{D+1}\;\frac{(-1)^{t+i-D}u_t(\theta_i)}{D+1}\;\frac{(-1)^{t+j-D}u_t(\theta_j)}{D+1},\nonumber
    \end{align}
    which simplifies to
    \begin{align}
        \frac{(-1)^{i+j+h}}{(D+1)^3}\sum_{t=0}^D(2t+1)u_t(\theta_h)u_t(\theta_i)u_t(\theta_j).\label{eq:BE-LHS-KT}
    \end{align} 
    Setting (\ref{eq:BE-LHS-KT}) equal to the right hand side of (\ref{eq:using-BE}) 
    and dividing by the coefficients completes the proof.
\end{proof}

\begin{corollary}\label{cor:phij-form-racah}
    For $0\leq h,i,j\leq D$ we have
    \begin{equation}\label{eq:phij-form-racah}
        p^h_{i,j}=(2i+1)(2j+1)(D+1)\Bigl(W\left({\textstyle \frac{D}{2},\frac{D}{2},i,h;j,\frac{D}{2}}\right)\Bigr)^2.
    \end{equation}
\end{corollary}

\begin{proof} Using Propositions \ref{prop:phij-ui-triple-product}, \ref{prop:biedenharn-elliott-kresch} and substituting (\ref{kidef}),(\ref{eq:nu}) we have
    \begin{align}
p^h_{i,j}&=\frac{k_ik_j}{\nu}\sum\limits_{t=0}^{D}k_t u_t(\theta_i)u_t(\theta_j)u_t(\theta_h)\nonumber\\
&=\frac{(2i+1)(2j+1)}{(D+1)^2}\paren{(D+1)^3\Bigl(W\left({\textstyle \frac{D}{2},\frac{D}{2},i,h;j,\frac{D}{2}}\right)\Bigr)^2}
\nonumber\\
&=(2i+1)(2j+1)(D+1)\Bigl(W\left({\textstyle \frac{D}{2},\frac{D}{2},i,h;j,\frac{D}{2}}\right)\Bigr)^2.
\nonumber
\end{align}

\end{proof}

\begin{corollary}\label{phijnonneg} For $0 \leq h,i,j \leq D$ we have $$p^h_{i,j} \geq 0.$$
\end{corollary}

\begin{proof}
Immediate from Corollary~\ref{cor:phij-form-racah}.
\end{proof}

\bigskip

\section{Proof of the Kresch-Tamvakis conjecture}  

    We are now ready to prove our main result. We will use the Perron-Frobenius theorem \cite[p.~529]{horn2012matrix}.

\begin{proposition}\label{prop:ui-bound}
     For $0 \leq i,j \leq D$ we have
$$|u_i(\theta_j)| \leq 1.$$
\end{proposition} 
\begin{proof}
By Lemma~\ref{lem:bi-eig-ki}, the vector ${\mathbb 1}$ is an eigenvector for $B_i$ with eigenvalue $k_i$. By Corollary \ref{phijnonneg}, the entries of $B_i$ are all nonnegative.
    By Lemma~\ref{lem:PBi-BisP} the scalar $v_i(\theta_j)$ is an eigenvalue of $B_i$. By the Perron-Frobenius theorem \cite[p.~529]{horn2012matrix}, 
    we have $|v_i(\theta_j)| \leq k_i$. The result follows from this and (\ref{uvrelnTH}). 
\end{proof}

\noindent
Equation (\ref{uijF}) and Proposition~\ref{prop:ui-bound} imply Theorem~\ref{thm:main}.

\bigskip 

\section{Appendix}

  In this appendix we give more detail about the formula for $p^h_{i,j}$ in Corollary~\ref{cor:phij-form-racah}. 
 By Lemma~\ref{lem:phij-sym}, without loss of generality we assume $i\leq j\leq h$. Also, in order to avoid trivialities we assume that $h,i,j$ satisfy the triangle inequalities; which in this case become $h\leq i+j$.
 As we evaluate $p^h_{i,j}$ in line (\ref{eq:phij-form-racah}) we consider the last factor. We evaluate that factor using Definition~\ref{Wdef} with

 $$a=\frac{D}{2},\qquad b=\frac{D}{2}, \qquad c=i, \qquad d=h, \qquad e=j,\qquad f=\frac{D}{2}.$$
 For these values,

        \begin{equation*}
         \alpha_1=D+i, \qquad \alpha_2=D+j, \qquad \alpha_3=D+h, \qquad\alpha_4=h+i+j,
     \end{equation*}
        
        \begin{equation*}
          \beta_1=D+i+j, \qquad \beta_2=D+h+i, \qquad \beta_3=D+h+j.
     \end{equation*}

\noindent
Note that
     \begin{equation*}
          \alpha_1-\beta_1=-j, \qquad \alpha_2-\beta_1=-i, \qquad \alpha_3-\beta_1=h-i-j, \qquad \alpha_4-\beta_1=h-D
     \end{equation*}
     
     \begin{equation*}
        -\beta_1-1=-D-i-j-1, \qquad \beta_2-\beta_1+1=h-j+1, \qquad \beta_3-\beta_1+1=h-i+1 .
     \end{equation*}

\medskip

\noindent
For the above data, (\ref{eq:phij-form-racah}) becomes

\begin{equation*}
    p^h_{i,j}=C^h_{i,j} (2i+1)(2j+1)(D+1)\,\Biggl({}_4F_3\left[\begin{matrix}
				-j ,\, -i ,\, h-i-j ,\, h-D\\
				-D-i-j-1 ,\, h-j+1 ,\, h-i+1 &
			\end{matrix};\,1\right]\Biggr)^2,
\end{equation*}
where 

\begin{align*}
    C^h_{i,j}&=\Bigg(\frac{\Delta({\textstyle\frac{D}{2},\frac{D}{2},i})\Delta({\textstyle\frac{D}{2},\frac{D}{2},j})\Delta({\textstyle\frac{D}{2},\frac{D}{2},h})\Delta(i,j,h))(D+i+j+1)!}{(h-i)!(h-j)!i!j!(i+j-h)!(D-h)!}\Bigg)^2\nonumber\\
    &=\frac{(D-i)!(D-j)!(D-h)!(j+h-i)!(h+i-j)!}{(D+i+1)!(D+j+1)!(D+h+1)!(i+j+h+1)!(i+j-h)!}\paren{\frac{h!(D+i+j+1)!}{(h-i)!(h-j)!(D-h)!}}^2.
\end{align*} 
 
\bigskip 

\noindent
{\bf Acknowledgement.} We would like to express our gratitude to Professor Paul Terwilliger, whose careful feedback greatly enhanced the clarity of the exposition.

\bibliographystyle{plain}
\bibliography{master}

\bigskip

\noindent
John S. Caughman\\
Fariborz Maseeh Dept of Mathematics \& Statistics\\
PO Box 751\\
Portland State University\\
Portland, OR 97207 USA\\
email: {\tt caughman@pdx.edu}

\bigskip

\noindent
Taiyo S. Terada\\
Fariborz Maseeh Dept of Mathematics \& Statistics\\
PO Box 751\\
Portland State University\\
Portland, OR 97207 USA\\
email: {\tt taiyo2@pdx.edu}

\end{document}